\magnification=1200

\hsize=11.25cm    
\vsize=18cm       
\parindent=12pt   \parskip=5pt     

\hoffset=.5cm   
\voffset=.8cm   

\pretolerance=500 \tolerance=1000  \brokenpenalty=5000

\catcode`\@=11

\font\eightrm=cmr8         \font\eighti=cmmi8
\font\eightsy=cmsy8        \font\eightbf=cmbx8
\font\eighttt=cmtt8        \font\eightit=cmti8
\font\eightsl=cmsl8        \font\sixrm=cmr6
\font\sixi=cmmi6           \font\sixsy=cmsy6
\font\sixbf=cmbx6

\font\tengoth=eufm10 
\font\eightgoth=eufm8  
\font\sevengoth=eufm7      
\font\sixgoth=eufm6        \font\fivegoth=eufm5

\skewchar\eighti='177 \skewchar\sixi='177
\skewchar\eightsy='60 \skewchar\sixsy='60

\newfam\gothfam           \newfam\bboardfam

\def\tenpoint{
  \textfont0=\tenrm \scriptfont0=\sevenrm \scriptscriptfont0=\fiverm
  \def\rm{\fam\z@\tenrm}
  \textfont1=\teni  \scriptfont1=\seveni  \scriptscriptfont1=\fivei
  \def\oldstyle{\fam\@ne\teni}\let\old=\oldstyle
  \textfont2=\tensy \scriptfont2=\sevensy \scriptscriptfont2=\fivesy
  \textfont\gothfam=\tengoth \scriptfont\gothfam=\sevengoth
  \scriptscriptfont\gothfam=\fivegoth
  \def\goth{\fam\gothfam\tengoth}
  
  \textfont\itfam=\tenit
  \def\it{\fam\itfam\tenit}
  \textfont\slfam=\tensl
  \def\sl{\fam\slfam\tensl}
  \textfont\bffam=\tenbf \scriptfont\bffam=\sevenbf
  \scriptscriptfont\bffam=\fivebf
  \def\bf{\fam\bffam\tenbf}
  \textfont\ttfam=\tentt
  \def\tt{\fam\ttfam\tentt}
  \abovedisplayskip=12pt plus 3pt minus 9pt
  \belowdisplayskip=\abovedisplayskip
  \abovedisplayshortskip=0pt plus 3pt
  \belowdisplayshortskip=4pt plus 3pt 
  \smallskipamount=3pt plus 1pt minus 1pt
  \medskipamount=6pt plus 2pt minus 2pt
  \bigskipamount=12pt plus 4pt minus 4pt
  \normalbaselineskip=12pt
  \setbox\strutbox=\hbox{\vrule height8.5pt depth3.5pt width0pt}
  \let\bigf@nt=\tenrm       \let\smallf@nt=\sevenrm
  \normalbaselines\rm}

\def\eightpoint{
  \textfont0=\eightrm \scriptfont0=\sixrm \scriptscriptfont0=\fiverm
  \def\rm{\fam\z@\eightrm}
  \textfont1=\eighti  \scriptfont1=\sixi  \scriptscriptfont1=\fivei
  \def\oldstyle{\fam\@ne\eighti}\let\old=\oldstyle
  \textfont2=\eightsy \scriptfont2=\sixsy \scriptscriptfont2=\fivesy
  \textfont\gothfam=\eightgoth \scriptfont\gothfam=\sixgoth
  \scriptscriptfont\gothfam=\fivegoth
  \def\goth{\fam\gothfam\eightgoth}
  
  \textfont\itfam=\eightit
  \def\it{\fam\itfam\eightit}
  \textfont\slfam=\eightsl
  \def\sl{\fam\slfam\eightsl}
  \textfont\bffam=\eightbf \scriptfont\bffam=\sixbf
  \scriptscriptfont\bffam=\fivebf
  \def\bf{\fam\bffam\eightbf}
  \textfont\ttfam=\eighttt
  \def\tt{\fam\ttfam\eighttt}
  \abovedisplayskip=9pt plus 3pt minus 9pt
  \belowdisplayskip=\abovedisplayskip
  \abovedisplayshortskip=0pt plus 3pt
  \belowdisplayshortskip=3pt plus 3pt 
  \smallskipamount=2pt plus 1pt minus 1pt
  \medskipamount=4pt plus 2pt minus 1pt
  \bigskipamount=9pt plus 3pt minus 3pt
  \normalbaselineskip=9pt
  \setbox\strutbox=\hbox{\vrule height7pt depth2pt width0pt}
  \let\bigf@nt=\eightrm     \let\smallf@nt=\sixrm
  \normalbaselines\rm}

\tenpoint

\def\pc#1{\bigf@nt#1\smallf@nt}         \def\pd#1 {{\pc#1} }

\catcode`\;=\active
\def;{\relax\ifhmode\ifdim\lastskip>\z@\unskip\fi
\kern\fontdimen2  -1.2 \fontdimen3 \string;}

\catcode`\:=\active
\def:{\relax\ifhmode\ifdim\lastskip>\z@\unskip\fi\penalty\@M\ \fi\string:}

\catcode`\!=\active
\def!{\relax\ifhmode\ifdim\lastskip>\z@
\unskip\fi\kern\fontdimen2  -1.1 \fontdimen3 \string!}

\catcode`\?=\active
\def?{\relax\ifhmode\ifdim\lastskip>\z@
\unskip\fi\kern\fontdimen2  -1.1 \fontdimen3 \string?}

\frenchspacing

\def\raggedbottom{\topskip 10pt plus 36pt\r@ggedbottomtrue}

\def\pointir{\unskip . --- \ignorespaces}

\def\Medbreak{\vskip-\lastskip\medbreak}

\long\def\th#1 #2\enonce#3\endth{
   \Medbreak\noindent
   {\pc#1} {#2\unskip}\pointir{\it #3}\smallskip}

\def\proof{\vskip-\lastskip\smallskip\noindent
 {\it Proof} : }

\def\decale#1{\smallbreak\hskip 28pt\llap{#1}\kern 5pt}
\def\decaledecale#1{\smallbreak\hskip 34pt\llap{#1}\kern 5pt}
\def\puce{\smallbreak\hskip 6pt{$\scriptstyle\bullet$}\kern 5pt}

\def\eqalign#1{\null\,\vcenter{\openup\jot\m@th\ialign{
\strut\hfil$\displaystyle{##}$&$\displaystyle{{}##}$\hfil
&&\quad\strut\hfil$\displaystyle{##}$&$\displaystyle{{}##}$\hfil
\crcr#1\crcr}}\,}

\catcode`\@=12

\showboxbreadth=-1  \showboxdepth=-1

\newcount\numerodesection \numerodesection=1
\def\section#1{\bigbreak
 {\bf\number\numerodesection.\ \ #1}\nobreak\medskip
 \advance\numerodesection by1}

\mathcode`A="7041 \mathcode`B="7042 \mathcode`C="7043 \mathcode`D="7044
\mathcode`E="7045 \mathcode`F="7046 \mathcode`G="7047 \mathcode`H="7048
\mathcode`I="7049 \mathcode`J="704A \mathcode`K="704B \mathcode`L="704C
\mathcode`M="704D \mathcode`N="704E \mathcode`O="704F \mathcode`P="7050
\mathcode`Q="7051 \mathcode`R="7052 \mathcode`S="7053 \mathcode`T="7054
\mathcode`U="7055 \mathcode`V="7056 \mathcode`W="7057 \mathcode`X="7058
\mathcode`Y="7059 \mathcode`Z="705A


\def\hfl#1#2#3{\smash{\mathop{\hbox to#3{\rightarrowfill}}\limits
^{\textstyle#1}_{\textstyle#2}}}

\def\agoth{{\goth a}}

\def\ogoth{{\goth o}}

\def\pgoth{{\goth p}}

\def\Q{{\bf Q}}

\def\Z{{\bf Z}}

\def\F{{\bf F}}

\def\Ker{\mathop{\rm Ker}\nolimits}

\def\normressym(#1,#2)_#3{\displaystyle\left({#1,#2\over#3}\right)}

\def\mod{\mathop{\rm mod.}\nolimits}
\def\pmod#1{\;(\mod#1)}

\newcount\refno 
\long\def\ref#1:#2<#3>{                                        
\global\advance\refno by1\par\noindent                              
\llap{[{\bf\number\refno}]\ }{#1} \pointir{\it #2} #3\goodbreak }

\def\citer#1(#2){[{\bf\number#1}\if#2\empty\relax\else,\ {#2}\fi]}

\def\pibar{{\bar\pi}}
\def\dii{d_2}

\newbox\bibbox
\setbox\bibbox\vbox{\bigbreak
\centerline{{\pc BIBLIOGRAPHIC} {\pc REFERENCES}}

\ref{\pc GAUSS} (C.):
Disquisitiones arithmeticae,
<Gerh.~Fleischer, Lipsiae, 1801, xviii+668 pp.>
\newcount\gauss \global\gauss=\refno


\ref{\pc HENSEL} (K.):
Die multiplikative Darstellung der algebraischen Zahlen f{\"u}r den Bereich
eines beliebigen Primteilers,
<J.\ f.\ d.\ reine und angewandte Math., {\bf 146}, 1916,
pp.~189--215.> 
\newcount\henselmult \global\henselmult=\refno


%

} 

\centerline{\bf Wilson's theorem} 
\bigskip\bigskip
\centerline{Chandan Singh Dalawat}
\bigskip

\centerline{\eightpoint{\it ...puisque de tels hommes n'ont pas cru ce
sujet indigne de leurs m{\'e}ditations...} \citer\gauss().\qquad}

\bigskip

More than two hundered years ago, Gauss generalised Wilson's theorem
($(p-1)!\equiv-1\pmod p$ for a prime number $p$) to an arbitrary integer $A>0$
in \S78 of his {\it Disquisitiones\/}~:

\th THEOREM 1 (\citer\gauss())
\enonce
Poductum ex omnibus numeris, numero\/ {\rm quocunque} dato $A$ minoribus
simulque ad ipsum primis, congruum est secundum $A$, vnitati vel negatiue vel
positiue sumtae.
\endth
(The product of all elements in $(\Z/A\Z)^\times$ is $\bar1$ or $-\bar1$).  He
then specifies that the product in question is $-\bar1$ if $A$ is $4$, or
$p^m$, or $2p^m$ for some odd prime $p$ and integer $m>0$~; it equals $\bar1$
in the remaining cases.

{\eightpoint According to Gauss (\citer\gauss(), \S76) the elegant theorem
according to which ``\thinspace upon augmenting the product of all numbers
less than a given prime number by the unity, it becomes divisible by that
prime number\thinspace'' was first stated by Waring in his {\it
Meditationes\/} --- which appeared in Cambridge in 1770 --- and attributed to
Wilson, but neither could prove it.  Waring remarks that the proof must be all
the more difficult as there is no {\it notation\/} which might express a prime
number.  {\it Nach unserer Meinung aber m{\"u}ssen derartige Wahrheiten
vielmehr aus Begriffen (}notionibus{\it ) denn aus Bezeichnungen
(}notationibus{\it ) gesch{\"o}pft werden} \citer\gauss().  The first proof
was given by Lagrange (1771).}


Some hundred years later, Hensel \citer\henselmult() developed his local
notions, which could have allowed him to extend the result from $\Z$ to rings
of integers in number fields~; our aim here is to show how he could have
done it.

\th PROPOSITION 2 (``Wilson's theorem'')
\enonce
For an ideal\/ $\agoth\subset\ogoth$ in the ring of integers of a number
field, the product of all elements in\/ $(\ogoth/\agoth)^\times$ is\/ $\bar1$,
except that it is

$-\bar1$ when\/ $\agoth$ has precisely one odd prime divisor, and
$v_\pgoth(\agoth)<2$ for every even prime ideal\/ $\pgoth$,

$\bar1+\pibar$ (resp.~$\bar1+\pibar^2$) when all prime divisors of\/ $\agoth$
are even and\/ for precisely one of them, say\/ $\pgoth$, $v_\pgoth(\agoth)>1$
with moreover\/ $v_\pgoth(\agoth)=2$, $f_\pgoth=1$
(resp.~$v_\pgoth(\agoth)=3$, $f_\pgoth=1$, $e_\pgoth>1$)~; here\/ $\pi$ is any
element of\/ $\pgoth$ not in\/ $\pgoth^2$.
\endth
The notation and the terminology are unambiguous~: a prime ideal $\pgoth$ of
$\ogoth$ is even if $2\in\pgoth$, odd if $2\notin\pgoth$~; $v_\pgoth(\agoth)$
is the exponent of $\pgoth$ in the prime decomposition of $\agoth$~;
$f_\pgoth$ is the residual degree and $e_\pgoth$ the ramification index of
$K_\pgoth|\Q_p$ ($p$ being the rational prime which belongs to $\pgoth$).

The proof will make it clear how to compute the product of all elements in any
given subgroup of $(\ogoth/\agoth)^\times$.

(It may happen that $\bar1+\pibar=-\bar1$ in $(\ogoth/\pgoth^2)^\times$
(resp.~$\bar1+\pibar^2=-\bar1$ in $(\ogoth/\pgoth^3)^\times$) for some even
prime $\pgoth\subset\ogoth$.  Example~: $\ogoth=\Z$ (resp.~$\Z[\sqrt2]$) and
$\pgoth$ the unique even prime of $\ogoth$.  More banally, we have $-\bar
1=\bar1$ in $(\ogoth/\pgoth^n)^\times$ when $\pgoth$ is an even prime and $n$
is between $1$ and $e_\pgoth$.)

\bigbreak 
\centerline{\bf 1. $d_2$} 
\medskip 

The elementary observation behind Gauss's proof of th.~1, also used in our
proof of prop.~2, is that the sum $s$ of all the elements in a finite
commutative group $G$ is $0$, unless $G$ has precisely one order-$2$ element
$\tau$, in which case $s=\tau$.  Anyone can supply a proof~; he can then skip
this section, and take the condition ``\thinspace$\dii(G)=1$\thinspace'' as a
shorthand for ``\thinspace$G$ has precisely one order-$2$ element\thinspace''.

Define $\dii(G)=\dim_{\F_2}({}_2G)$, where ${}_2G$ is the subgroup of $G$
killed by~$2$.  It is clear that $G$ has $2^{\dii(G)}-1$ order-$2$ elements.

\th EXAMPLE 3
\enonce
For a prime number\/ $p$ and integer\/ $n>0$, we have
$$
\dii((\Z/p^n\Z)^\times)=
\cases{
1&if\/ $p\neq2$,\cr
0&if\/ $p=2$ and $n=1$,\cr
1&if\/ $p=2$ and $n=2$,\cr
2&if\/ $p=2$ and $n>2$.\cr
}
$$
\endth

In this example, the unique order-$2$ element is $-\bar1$  whenever $\dii=1$.

\th LEMMA 4
\enonce
The sum $s$ of all elements in the group\/ $G$ is $0$ unless $\dii(G)=1$, in
which case $s$ is the unique order-$2$ element of\/ $G$.
\endth
The involution $\iota:g\mapsto-g$ fixes every element of the subgroup
${}_2G=\Ker(x\mapsto2x)$.  As the sum of elements in the remaining orbits of
$\iota$ is~$0$, we are reduced to the case $G={}_2G$ of a vector $\F_2$-space,
and the proof is over by induction on the dimension $\dii(G)$ of ${}_2G$,
starting with dimension~$2$.

\medbreak

{\it Proof of Gauss's th.~1}~: Let $A=\prod_pp^{m_p}$ be the prime
decomposition of $A$.  By the Chinese remainder theorem, $(\Z/A\Z)^\times$ is
the product over $p$ of $(\Z/p^{m_p}\Z)^\times$, so $\dii((\Z/A\Z)^\times)$ is
the sum over $p$ of $\dii((\Z/p^{m_p}\Z)^\times)$.  In view of example~3, the
only way for this sum to be $1$ is for $A$ to be $2^2$, or $p^{m_p}$, or
$2p^{m_p}$ for some odd prime $p$ and integer $m_p>0$.

\bigbreak 
\centerline{\bf 2. Local units} 
\medskip 

Let's enter Hensel's world~: let\/ $p$ be a prime number, $K\,|\,\Q_p$ a
finite extension, $\ogoth$ its ring of integers, $\pgoth$ the unique maximal
ideal of\/ $\ogoth$.  Let\/ $n>0$ be an integer.  We would like to know when
$\dii((\ogoth/\pgoth^n)^\times)=1$, and, when such is the case, which one the
unique order-$2$ element is.

\th PROPOSTION 5
\enonce
Denoting by\/ $e$ the ramification index and by\/ $f$ the residual degree of\/ 
$K\,|\,\Q_p$,  
$$\def\\{\;\phantom>}
\dii((\ogoth/\pgoth^n)^\times)=
\cases{\\1&if\/ $p\neq2$,\cr
\\0&if\/ $p=2$, $n=1$,\cr
\\1&if\/ $p=2$, $n=2$, $f=1$,\cr
\\1&if\/ $p=2$, $n=3$, $f=1$, $e>1$,\cr
>1&in all other cases.\cr}
$$
For any\/ $\ogoth$-basis\/ $\pi$ of\/ $\pgoth$, the unique order-$2$ element
in the cases\/ $\dii=1$ is 
$$
\cases{
-\bar1&if\/ $p\neq2$,\cr
\bar1+\pibar&if\/ $p=2$, $n=2$, $f=1$,\cr
\bar1+\pibar^2&if\/ $p=2$, $n=3$, $f=1$, $e>1$.\cr
}
$$
\endth
\proof For every $j>0$, denote by $U_j$ the kernel of
$\ogoth^\times\rightarrow(\ogoth/\pgoth^j)^\times$.  If $p\neq2$, the group
$(\ogoth/\pgoth^n)^\times$ is the direct product of the even-order cyclic
group $(\ogoth/\pgoth)^\times$ and the $p$-group $U_1/U_n$, so $\dii=1$.

Assume now that $p=2$.  When $n=1$, the group $(\ogoth/\pgoth)^\times$ is
(cyclic) of odd order, so $\dii=0$.  If $f>1$, then the $\dii$ of $U_1/U_2$ is
$f$ and hence the $\dii$ of $(\ogoth/\pgoth^n)^\times$ is $>1$ for every
$n>1$.

Assume further that $f=1$.  When $n=2$, the $\dii$ of
$(\ogoth/\pgoth^2)^\times=U_1/U_2$ is $f=1$.  If moreover $e=1$, then the
$\dii$ of $U_1/U_n$ is $2$ for $n>2$ (example~3).

Assume finally that, in addition, $e>1$.  We see that $U_1/U_3$ is generated
by $\bar1+\pibar$, since $(1+\pi)^2=1+\pi^2+2\pi$ is in
$U_2$ but not in $U_3$. However, $U_1/U_4$ is not cyclic because its order is
$8$ whereas every element has order at most~$4$~: for every $a\in\ogoth$,
$$
(\bar1+\bar a\pibar)^4=
\bar1+\bar4\pibar\bar a+\bar6\pibar^2\bar a^2
+\bar4\pibar^3\bar a^3+\pibar^4\bar a^4=\bar1
$$ 
in $U_1/U_4$.  Hence $U_1/U_n$ is not cyclic for $n>3$ (cf.~Narkiewicz, {\it
Elem.\ and anal.\ theory of alg.\ numbers}, 1990, p.~275). This concludes the
proof.

(For $p=2$ and $n>2e$, we have $\dii((\ogoth/\pgoth^n)^\times)=1+ef$~;
cf.~Hasse, {\it Zahlentheorie}, Kap.~15.)

\th COROLLARY 6
\enonce
The only cases in which the group\/ $(\ogoth/\pgoth^n)^\times$ has precisely
one order-$2$ element $s$ are~: $p\neq2$~; $p=2$, $n=2$, $f=1$~; $p=2$, $n=3$,
$f=1$, $e>1$.  In these three cases, $s=-\bar1$, $\bar1+\pibar$,
$\bar1+\pibar^2$, respectively.  The group\/ $(\ogoth/\pgoth^n)^\times$ has no
order-$2$ element precisely when $p=2$, $n=1$.
\endth

\bigbreak 
\centerline{\bf 2. The proof} 
\medskip 

Let us return to the global situation of an ideal $\agoth\subset\ogoth$ in the
ring of integers of a number field $K\,|\,\Q$.  The proof can now proceed as
in the case $\ogoth=\Z$ (\S1).  Everything boils down to deciding if the
$\dii$ of $(\ogoth/\agoth)^\times$ is $1$ --- we know that the product of all
elements is $1$ if $\dii\neq1$ (lemma~4).  Writing
$\agoth=\prod_\pgoth\pgoth^{m_\pgoth}$ the prime decomposition of $\agoth$,
the Chinese remainder theorem tells us that $\dii((\ogoth/\agoth)^\times)$ is
the sum, over the various primes $\pgoth$ of $\ogoth$, of
$\dii((\ogoth/\pgoth^{m_\pgoth})^\times)$.  This sum can be $1$ only when one
of the terms is $1$, the others being $0$.

For each $\pgoth$, the group $(\ogoth/\pgoth^{m_\pgoth})^\times$ is the same
as $(\ogoth_\pgoth/\pgoth_\pgoth^{m_\pgoth})^\times$, where $\ogoth_\pgoth$ is
the completion of $\ogoth$ at $\pgoth$ and $\pgoth_\pgoth$ is the unique
maximal ideal of $\ogoth_\pgoth$.  Running through the possibilities
enumerarted in prop.~5 completes the proof of prop.~2.

{\it Example\/~7}.  Let $\zeta\in\bar\Q^\times$ be an element of order~$2^t$
($t>1$)~; take $K=\Q(\zeta)$ and $\pgoth$ the unique even prime of its ring of
integers $\Z[\zeta]$.  We have $e_\pgoth=2^{t-1}$ and $f_\pgoth=1$~; we may
take $\pi=1-\zeta$.  The product of all elements in
$(\Z[\zeta]/\pgoth^n)^\times$ is respectively $\bar1$, $\bar1+\pibar$,
$\bar1+\pibar^2$, $\bar1$ for $n=1$, $n=2$, $n=3$ and $n>3$.

\bigbreak 
\centerline{\bf 3. Acknowledgements} 
\medskip 

{\eightpoint We thank Herr Prof.~Dr.~Peter Roquette for suggesting the present
definition $\dii(G)=\dim_{\F_2}({}_2G)$ instead of the original
$\dii(G)=\dim_{\F_2}(G/2G)$.  After this Note was completed, a search in the
literature revealed M.~La{\v s}{\v s}{\'a}k, {\it Wilson's theorem in
algebraic number fields}, Math.\ Slovaca, {\bf 50} (2000), no.~3,
pp.~303--314.  We solicited a copy from Prof.~G.~Grekos, and thank him for
supplying one~; it contains substantially the same result as our prop.~2.  
Our proof is shorter, simpler, more direct, and more conceptual~; it is based
on {\it notionibus\/} rather than {\it notationibus}, of which there is
now-a-days a surfeit.  In any case, our aim was to show how Hensel could have
proved prop.~2. }


\bigbreak
\unvbox\bibbox 
\vskip2cm plus 1cm minus 1cm
\vbox{{\obeylines\parskip=0pt\parindent=0pt
Chandan Singh Dalawat
Harish-Chandra Research Institute
Chhatnag Road, Jhunsi
{\pc ALLAHABAD} 211\thinspace019, India
\vskip5pt
\tt dalawat@gmail.com}}

\bye